\documentclass[reqno]{amsart}
\usepackage{amsmath,mathtools,amssymb}

\mathtoolsset{showonlyrefs}

\usepackage{hyperref}
\usepackage[normalem]{ulem}

\newtheorem{question}{Question}

\theoremstyle{definition}

\newcommand{\R}{\mathbb{R}}
\newcommand{\C}{\mathbb{C}}

\newcommand\re{\operatorname{Re}}
\newcommand\im{\operatorname{Im}}

\newcommand\Tr{\operatorname{Tr}}

\newcommand\eps\varepsilon
\renewcommand\epsilon\varepsilon

\newcommand{\dist}{\operatorname{dist}}
\newcommand{\rd}{\mathrm{d}}

\title[]{Open problem: Violation of locality for Schr\"odinger operators with complex potentials}

\author{Jean-Claude Cuenin}
\address[Jean-Claude Cuenin]{Department of Mathematical Sciences, Loughborough University, Loughborough, Leicestershire, LE11 3TU United Kingdom}
\email{J.Cuenin@lboro.ac.uk}

\author{Rupert L. Frank}
\address[Rupert L. Frank]{Mathe\-matisches Institut, Ludwig-Maximilians Universit\"at M\"unchen, The\-resienstr.~39, 80333 M\"unchen, Germany, and Munich Center for Quantum Science and Technology, Schel\-ling\-str.~4, 80799 M\"unchen, Germany, and Mathematics 253-37, Caltech, Pasa\-de\-na, CA 91125, USA}
\email{r.frank@lmu.de}

\date{September 17, 2024}

\thanks{\copyright\, 2024 by the authors. This paper may be reproduced, in its entirety, for non-commercial purposes.\\
	Partial support through US National Science Foundation (DMS-1954995; RLF), as well as through the German Research Foundation (EXC-2111-390814868 and TRR 352-470903074; RLF) is acknowledged. Support through the Engineering \& Physical Sciences Research Council (EP/X011488/1; JCC) is acknowledged.}

\begin{document}

\begin{abstract}
We explain in which sense Schr\"odinger operators with complex potentials appear to violate locality (or Weyl's asymptotics), and we pose three open problems related to this phenomenon.
\end{abstract}
\maketitle

\section{Background}

We are interested in the semiclassical Schr\"odinger operator
\begin{align}
	-\hbar^2\Delta+V\quad \mbox{in}\quad L^2(\R^d)
\end{align}
with \emph{complex-valued} potential $V$ and semiclassical parameter $\hbar\in (0,1]$. We will assume throughout that the potential $V$ is locally sufficiently regular, so that the operator can be defined as an $m$-sectorial operator, and that it decays at infinity (at least in some averaged sense), so that the spectrum of the operator in $\C \setminus [0, \infty)$ consists of isolated eigenvalues of finite algebraic multiplicities. We denote these eigenvalues by $E_j$, repeating an eigenvalue according to its algebraic multiplicity. We seek bounds on sums of the $E_j$ in terms of integrals of the potential $V$.

\subsection{Real-valued potentials}
We begin by reviewing the case of real-valued potentials $V$. If $V$ is, say, continuous and compactly supported, then Weyl's asymptotic formula says that, for $\gamma\geq 0$,
\begin{align}\label{semiclassical Riesz means real V}
    \lim_{\hbar\to 0}\hbar^d\Tr(-\hbar^2\Delta+V)_-^{\gamma}=L_{\gamma,d}^{\rm cl}\int_{\R^d}V(x)_-^{\gamma+d/2}\rd x,
\end{align}
where $L_{\gamma,d}^{\rm cl}=\int_{\R^d}(|\xi|^2-1)_-^{\gamma}\frac{\rd\xi}{(2\pi)^d}$ is the semiclassical constant (see \cite{MR4496335}) and $V(x)_-=\max(0,-V(x))$ is the negative part of $V(x)$. For $\gamma=0$, the left-hand side of \eqref{semiclassical Riesz means real V} is interpreted as the number of negative eigenvalues of $-\hbar^2\Delta+V$, and for $\gamma>0$ it is called the Riesz mean of order $\gamma$.

A non-asymptotic bound that captures the correct order of magnitude of the Riesz means \eqref{semiclassical Riesz means real V} in the asymptotic regime would be of the form
\begin{align}\label{Lieb-Thirring real V}
\Tr(-\hbar^2 \Delta+V)_-^{\gamma} =  \sum_j|E_j|^{\gamma}\leq L_{\gamma,d} \hbar^{-d} \int_{\R^d}V(x)_-^{\gamma+d/2}\rd x,   
\end{align}
where $L_{\gamma,d}$ is some positive constant depending on $\gamma,d$, but independent of $V$ and $\hbar>0$. Note that necessarily $L_{\gamma,d}\geq L_{\gamma,d}^{\rm cl}$. For the sake of the following discussion, we have stated the inequality \eqref{Lieb-Thirring real V} for the semiclassical Schr\"odinger operator $-\hbar^2\Delta+V$, but a simple scaling argument shows that it is equivalent to the corresponding bound for $\hbar=1$.

The validity of the bound \eqref{Lieb-Thirring real V} for $\gamma>1/2$ for $d=1$ and for $\gamma>0$ in $d\geq 2$ is a celebrated result of Lieb and Thirring. The case $\gamma=0$ in $d\geq 3$ is due Cwikel, Lieb and Rozenblum, and the case $\gamma=1/2$ in $d=1$ is due to Weidl. We refer to \cite{MR4496335} for more background and a discussion of optimal constants. (For complex-valued potentials, we will have nothing to say about constants.)

\subsection{Complex-valued potentials}
We now turn our attention to complex-valued potentials. We are interested in bounds of a similar nature as \eqref{Lieb-Thirring real V}. However, it is not at all obvious what the most natural analogues of these bounds should look like. 
\begin{itemize}
    \item The naive generalization
\begin{align}\label{Lieb-Thirring complex V naive}
  \sum_j|E_j|^{\gamma}\leq C_{\gamma,d} \hbar^{-d} \int_{\R^d}|V(x)|^{\gamma+d/2}\rd x   
\end{align}
fails for $\gamma>1/2$, even for a single eigenvalue \cite{MR4561804}. 
\item Frank, Laptev, Lieb and Seiringer \cite{MR2260376} proved that \eqref{Lieb-Thirring complex V naive} is valid for all eigenvalues outside an arbitrary fixed cone with a constant that becomes unbounded as the angle of the cone tends to zero. More precisely, they proved that for $\gamma\geq 1$, $\kappa>0$,
\begin{align}\label{Frank--Laptev--Lieb--Seiringer}
  \sum_{|\im E_j|\geq \kappa \re E_j}|E_j|^{\gamma}\leq C_{\gamma,d}(1+\kappa^{-1})^{\gamma+d/2} \hbar^{-d} \int_{\R^d}|V(x)|^{\gamma+d/2}\rd x.   
\end{align}
B\"ogli \cite{MR4651274} has shown that for $d=1$ as $\kappa\to 0$ the order of divergence $\kappa^{-\gamma-d/2}$ of the constant in \eqref{Frank--Laptev--Lieb--Seiringer} is optimal.
\item Averaging the bound \eqref{Frank--Laptev--Lieb--Seiringer} with respect to the opening angle $\kappa$, Demuth--Hansmann--Katriel \cite{MR2559715} established the following inequality, valid for \emph{all} eigenvalues,
\begin{align}\label{DKH bound}
 \sum_j |E_j|^{-\sigma}\delta(E_j)^{\gamma+\sigma}\leq C_{\gamma,\sigma,d} \hbar^{-d} \int_{\R^d}|V(x)|^{\gamma+d/2}\rd x,
\end{align}
where $\gamma\geq 1$, $\sigma>d/2$ and where we set
$$
\delta(E):=\dist(E,[0,\infty))
$$
Since $\delta(E_j)/|E_j|\leq 1$, the bound \eqref{DKH bound} becomes stronger the smaller $\sigma$. Demuth--Hansmann--Katriel \cite{MR3016473} asked whether \eqref{DKH bound} remains true for $\sigma=d/2$ and $\gamma>0$ (with $\gamma\geq 1/2$ if $d=1$). B\"ogli and Štampach \cite{MR4322041} answered this question in the negative for $d=1$ and $\gamma\geq 1/2$.
\end{itemize}
We note that the inequalities \eqref{Frank--Laptev--Lieb--Seiringer} and \eqref{DKH bound} reduce to the standard Lieb--Thirring inequalities when $V$ is real-valued since in this case, $E_j<0$ and thus $\delta(E_j)=|E_j|$.

Our first question concerns the assumption $\gamma\geq 1$ in \eqref{DKH bound}.

\begin{question}
	Let $d\geq 1$, $0<\gamma<1$ (with $\gamma\geq1/2$ if $d=1$) and $\sigma>d/2$. Does there exist a constant $C_{\gamma,\sigma,d}$ such that for all $V\in L^{\gamma+d/2}(\R^d)$ the inequality
	\begin{align}
		\sum_j |E_j|^{-\sigma}\delta(E_j)^{\gamma+\sigma}\leq C_{\gamma,\sigma,d}\hbar^{-d}\int_{\R^d}|V(x)|^{\gamma+d/2}\rd x
	\end{align}   
	holds for all $\hbar>0$?
\end{question}

Note that in this question one may assume without loss of generality that $\hbar=1$.

\section{Violation of locality}

For real-valued potentials, a key feature of the asymptotics \eqref{semiclassical Riesz means real V}, which is captured by the Lieb-Thirring inequalities \eqref{Lieb-Thirring real V} and is one of the reasons for their usefulness in applications, is \emph{locality}: sums of eigenvalues are estimated by an integral involving the potential. Hence, two disjoint pieces of $V$ contribute additively to the asymptotics \eqref{semiclassical Riesz means real V} and to the upper bound \eqref{Lieb-Thirring real V}. This is closely related to the appearance of $\hbar^{-d}$ on the right side of the inequality. For complex-valued potentials, these features are retained by the inequalities \eqref{Frank--Laptev--Lieb--Seiringer} and \eqref{DKH bound}.

In this section we first discuss some examples in one dimension where eigenvalue sums are $\gg \hbar^{-d}$ as $\hbar\to 0$ and then we present Lieb--Thirring-type bounds of a different nature than \eqref{Frank--Laptev--Lieb--Seiringer} and \eqref{DKH bound}, which, to some extent, capture a growth that is faster than $\hbar^{-d}$. In the following discussion, we will focus on the power of $\hbar^{-1}$ as an indication of the validity or the degree of violation of locality.

\subsection{Examples of eigenvalue sums that are $\gg \hbar^{-d}$}

We consider eigenvalue sums as on the left side of \eqref{DKH bound}, that is, of the form $\sum_j |E_j|^{-\sigma}\delta(E_j)^\beta$, and recall some examples from the literature where they are not bounded from above by a constant times $\hbar^{-d}$. So far, these examples are limited to dimension $d=1$. 

It is shown in \cite[Theorem 4]{MR4426735} that there exists $V\in L^1(\R)$ and a constant $c>0$ such that for all sufficiently small $\hbar>0$,
\begin{align}\label{Cuenin example 1d}
	\sum_j \delta(E_j)\geq c \Big(\frac{\hbar^{-1}}{\log \frac{1}{\hbar}} \Big)^{2}.   
\end{align}
In the same vein, the example in \cite{MR4322041} shows that there exists $V\in L^1\cap L^\infty(\R)$ (in fact, $i$ times the characteristic function of an interval) such that for every $\beta\geq 1$ there is a constant $c_\beta>0$ such that for all sufficiently small $\hbar>0$,
\begin{align}\label{Bogli-Stampach}
	\sum_j |E_j|^{-1/2}\delta(E_j)^\beta \geq c_\beta \hbar^{-1} \log\frac{1}{\hbar}.
\end{align}
In fact, by looking at the computations in \cite{MR4322041} we see that for any $\beta\geq 1$, any $0\leq\sigma<1/2$ and any $\epsilon>0$ there is a constant $c_{\beta,\sigma,\epsilon}>0$ such that
\begin{align}\label{Bogli-Stampach2}
	\sum_j |E_j|^{-\sigma}\delta(E_j)^{\beta} \geq c_{\beta,\sigma,\epsilon} \hbar^{-2+2\sigma+\epsilon}.
\end{align}
For $\sigma=0$ and $\beta=1$ this is similar, but slightly worse than \eqref{Cuenin example 1d}.

In particular, for any $\gamma\geq 1/2$ there is a potential $V\in L^{\gamma+1/2}(\R)$ such that sums $\sum_j |E_j|^{-\sigma}\delta(E_j)^{\gamma+\sigma}$ are $\gg h^{-1}$ provided $\sigma\leq1/2$. For $\gamma\geq 1$ this should be contrasted with the bound \eqref{DKH bound}, which implies that this sum is $\lesssim \hbar^{-1}$ for any $\sigma>1/2$. In this sense, $\sigma=1/2$ appears to be a threshold for $\gamma\geq 1$. If the answer to Question 1 is affirmative in dimension $d=1$, it is also a threshold for $1/2\leq\gamma<1$.

It is natural to wonder whether the above examples have analogues in higher dimensions.

\begin{question}
	Let $d\geq 2$, $\gamma>0$ and $0\leq\sigma\leq d/2$. Does there exist a $V\in L^{\gamma+d/2}(\R^d)$ such that
		\begin{align}
			\sum_j |E_j|^{-\sigma}\delta(E_j)^{\gamma+\sigma}\gg \hbar^{-d}
		\end{align}   
		as $\hbar\to 0$?
\end{question}

Sabine B\"ogli has informed us that she has made progress towards an answer to this question.

One can also ask (for any $d\geq 1$) whether the violation of the $\hbar^{-d}$-bound is power-like (as in \eqref{Cuenin example 1d} and \eqref{Bogli-Stampach2}) or logarithmic (as in \eqref{Bogli-Stampach}). In the first case, one might ask for the optimal power. We will see below that the lower bound in \eqref{Cuenin example 1d} is optimal up to powers of logarithms.

One may wonder whether the violation of locality holds for any generic choice of $V$ or only for certain specially chosen ones. Similarly, to which extent does the power of $\hbar$ depend on the choice of the potential? Also, it would be interesting to understand the role of the purported threshold $\sigma=d/2$ for the validity/violation of an $\hbar^{-d}$-bound on $\sum_j |E_j|^{-\sigma}\delta(E_j)^{\gamma+\sigma}$.

An even vaguer question is to understand from a more conceptual point of view the deeper reason behind a potential violation of locality.

\subsection{Nonlocal Lieb--Thirring bounds}

A family of eigenvalue bounds that is different in nature from \eqref{Frank--Laptev--Lieb--Seiringer} and \eqref{DKH bound} has been obtained by Frank and Sabin \cite{MR3730931}, as well as by Frank \cite{MR3717979}. These bounds only retain the scale-invariance but lose locality. They are of the form
\begin{align}\label{Frank--Sabin type bounds}
	\sum_j |E_j|^{\alpha}\left(\frac{\delta(E_j)}{|E_j|}\right)^{\beta}\leq C_{\alpha,\beta,\gamma,d} \Big(\hbar^{-d}\int_{\R^d}|V(x)|^{\gamma+d/2}\rd x\Big)^{\alpha/\gamma}
\end{align}
for certain values of $\alpha,\beta,\gamma,d$ (see below). As before, the result is equivalent to the corresponding result for $\hbar =1$.

A general observation is that all known bounds of the type \eqref{Frank--Sabin type bounds} have $\alpha/\gamma>1$. This means that a sum of eigenvalues is bounded by a \emph{power} of an integral, and this power is strictly greater than $1$. This is what we mean by a \emph{loss of locality}. Disjoint pieces of the potential no longer contribute additively to the bound on an eigenvalue sum. It also means that the power of $\hbar^{-1}$ that appears in the bound on the eigenvalue sum is strictly larger than the semiclassical power $d$.


An example from \cite{MR3717979} is the following bound in $d=1$,
\begin{align}\label{Frank III d=1}
\sum_j \delta(E_j)\leq C\Big(\hbar^{-1}\int_{\R}|V(x)|\rd x\Big)^{2},
\end{align}
Comparing this upper bound with the lower bound in \eqref{Cuenin example 1d}, we see that \eqref{Frank III d=1} is sharp as $\hbar\to 0$ up to logarithms.

The following instances of \eqref{Frank--Sabin type bounds} have been proved by Frank--Sabin \cite{MR3730931} and Frank \cite{MR3717979}. Notice that in (c), (d), (e), (f), the sum is restricted to $j$ satisfying either $|E_j|^\gamma \leq \hbar^{-d} \int_{\R^d}|V|^{\gamma+d/2}\rd x$ or  $|E_j|\geq \hbar^{-d} \int_{\R^d}|V|^{\gamma+d/2}\rd x$.

\begin{itemize}
    \item[(a)] $\alpha=1/2$, $\beta=1$, $0<\gamma<d/(2(2d-1))$, $d\geq 2$ (\cite[Thereom 16]{MR3730931}, see also \cite[(1.5)]{MR3717979}).
    \item[(b)] $\alpha>(d-1)\gamma/(d/2-\gamma)$, $\beta=1$, $d/(2(2d-1))\leq \gamma\leq 1/2$, $d\geq 2$ (\cite[Thereom 16]{MR3730931}, see also \cite[(1.5)]{MR3717979}).
    \item[(c)] $\alpha=\beta>2\gamma$, $\gamma>1/2$, truncation $|E_j|^\gamma \leq \hbar^{-d} \int_{\R^d}|V|^{\gamma+d/2}\rd x$ \cite[(1.7)]{MR3717979}).
    \item[(d)] $0<\alpha<\gamma(\gamma+d/2)$, $\beta>2\gamma$, $\gamma>1/2$, truncation $|E_j|^\gamma \geq \hbar^{-d} \int_{\R^d}|V|^{\gamma+d/2}\rd x$ \cite[(1.8)]{MR3717979}). 
    \item[(e)] $\alpha=\beta=\gamma+d/2$, $\gamma\geq 1/2$ ($d=1)$ or $\gamma>0$ ($d\geq 2$), truncation $|E_j|^\gamma \leq \hbar^{-d} \int_{\R^d}|V|^{\gamma+d/2}\rd x$ \cite[(1.9)]{MR3717979}).
    \item[(f)] $\alpha>\gamma$, $\beta=\gamma+d/2$, $\gamma\geq 1/2$ ($d=1)$ or $\gamma>0$ ($d\geq 2$), truncation $|E_j|^\gamma \geq \hbar^{-d} \int_{\R^d}|V|^{\gamma+d/2}\rd x$ \cite[(1.10)]{MR3717979}).
\end{itemize}
Note that \eqref{Frank III d=1} is a particular case of (e) in $d=1$ since $|E_j|^{1/2} \leq C \hbar^{-1} \int_{\R}|V(x)|\rd x$ by \cite{MR1819914}. As we have seen, this inequality is sharp up to logarithms.

\begin{question}
 For which values of $d\geq 1$ and $\alpha,\beta,\gamma$ listed in (a)--(f) is \eqref{Frank--Sabin type bounds} sharp up to a factor of $\hbar^{\eps}$, for arbitrary $\eps>0$? Can one increase the parameter region where bounds of the type \eqref{Frank--Sabin type bounds} are valid?
\end{question}

\bibliographystyle{plain}
\bibliography{Bibliography.bib}
\end{document}